\documentclass[11pt,a4paper,twoside]{article}

\usepackage{diagrams,a4,indentfirst,theorem,latexsym,QED}

\date{16.08.2004}

\diagramstyle[labelstyle=\scriptstyle]
\newarrow {Cof}      {triangle}--->

\newcommand\nothing[1]{\relax}

\renewcommand{\qedsymbol}{\quad\Box}
\renewcommand{\TheWordProof}{\bf Proof.}

\theoremstyle{change}
\theorembodyfont{\it}
\newtheorem{theorem}{Theorem.}[subsection]
\newtheorem{proposition}[theorem]{Proposition.}
\newtheorem{lemma}[theorem]{Lemma.}
\newtheorem{corollary}[theorem]{Corollary.}
{\theorembodyfont{\upshape} \newtheorem{definition}[theorem]{Definition.}
 }

\newcommand{\ie}{{\it i.e.}}

\newcommand{\simp}{\hbox{\rm Simp\,}}
\newcommand{\simpg}{\hbox{\rm Simp}_G}
\newcommand{\bZ}{\mathbf{Z}}
\newcommand{\op}{\mathrm{op}}
\newcommand{\fun}{\mathrm{Fun}\,}
\newcommand{\ab}{\mathrm{ab}}

\begin{document}

\author{Thomas~H\"uttemann}

\title{The Linearisation Map in Algebraic $K$-Theory}
\maketitle
\centerline {\it Mathematisches Institut, Universit\"at G\"ottingen}
\centerline {\it Bunsenstr.~3--5, D--37073 G\"ottingen, Germany}
\centerline {e-mail: \texttt{huette@uni-math.gwdg.de}}

\vglue 1\bigskipamount \hrule \medskip

{\footnotesize\noindent Let $X$ be a pointed connected simplicial set with
  loop group~$G$. The linearisation map in $K$-theory as defined by
  \textsc{Waldhausen\/} uses $G$-equivariant spaces. This
  paper gives an alternative description using presheaves of sets and
  \textsc{abel}ian groups on the simplex category of~$X$. In other words, the
  linearisation map is defined in terms of~$X$ only, avoiding the use of the
  less geometric loop group. The paper also includes a comparison of
  categorical finiteness with the more geometric notion of finite $CW$ objects
  in cofibrantly generated model categories. The application to the
  linearisation map employs a model structure on the category of
  \textsc{abel}ian group objects of retractive spaces over~$X$.  \smallskip

\noindent
{\it AMS subject classification (2000):\/} primary 19D10, secondary 55U35

\noindent
{\it Keywords:\/} Algebraic $K$-Theory of spaces, linearisation, finiteness
conditions}\medskip \hrule \vglue 1\bigskipamount

\centerline{to appear in Forum Mathematicum}

\subsection*{Introduction}

Let $W$ be a simplicial set with a (right) action of the simplicial
monoid~$G$. The category of $G$-equivariant retractive spaces over~$W$ is
denoted $\mathcal{R}(W,G)$; objects are triples $(Y,r,s)$ where $Y$ is a
simplicial set equipped with a $G$-action, $r \colon Y \rTo W$ is a
$G$-equivariant map, and $s$~is a $G$-equivariant section of~$r$. If $G$ is
the trivial monoid, it is omitted from the notation.

The category $\mathcal{R}(W,G)$ has been introduced by \textsc{Waldhausen\/}
\cite{W-spaces} to study the algebraic $K$-theory of spaces. For a pointed
connected simplicial set~$X$, the $K$-theory space $A(X)$ is defined as
$A(X):= \Omega |h\mathcal{S}_\bullet \mathcal{R}(X)_\mathrm{f}| $
where $ \mathcal{R}(X)_\mathrm{f}$ is a certain subcategory of
$\mathcal{R}(X)$ consisting of ``finite'' objects. It is proved in
\cite{W-spaces} that one can replace $\mathcal{R}(X)$ by the category
$\mathcal{R}(*,G)$ where $G$ is a loop group of~$X$. The latter category is
the domain of the linearisation functor which defines a natural map from
$A(X)$ to the algebraic $K$-theory space $K(\bZ[G])$ of the simplicial group
ring $\bZ[G]$.

Since the loop group $G$ is a less geometric object than the simplicial
set~$X$ itself, it is useful to have a description of the linearisation map
purely in terms of spaces over~$X$ as suggested by \textsc{Waldhausen}
\cite[\S2.3, p.~400]{W-spaces}. The domain of the new linearisation map is the
category $\mathcal{R}(X)$, the codomain is the category $\mathcal{R}^\ab (X)$
of \textsc{abel}ian group objects in~$\mathcal{R}(X)$, and the map is induced
by applying the (reduced) free \textsc{abel}ian group functor to the fibres of
the retractions $Y \rTo X$.

The goal of this paper is to provide a detailed account of the linearisation
map. In particular, we discuss finiteness conditions for objects in
$\mathcal{R}^\ab (X)$ (not treated in \cite[\S2.3]{W-spaces}), and give an
intrinsic definition of weak equivalences in~$\mathcal{R}^\ab (X)$ (avoiding
the use of simplicial $\bZ[G]$-modules as in~{\it loc.cit.}).

\medbreak

In \S1 we introduce the relevant categories and functors. In particular, we
prove that retractive objects are nothing but presheaves on the simplex
category:

\medbreak
\noindent {\bf Proposition~\ref{FunIsRofX} and~\ref{FunIsAb}.}
{\it Suppose $W$ is a simplicial set equipped with an action of the simplicial
  group~$G$. Then there are equivalences of categories
  $\mathcal{R} (W,G) \cong \fun (\simpg (W)^\op,\, \mathrm {\rm Set}_* )$
  and $\mathcal{R}^\ab (W,G) \cong \fun (\simpg (W)^\op,\, \mathrm
  {\rm Ab})$.
  Here $\simpg(W)$ is the equivariant simplex category of~$W$
  (Definition~\ref{def:simplexcat}), and the superscript ``$\ab$'' denotes
  \textsc{abel}ian group objects.}
\medbreak

In \S2 we set up the necessary machinery to discuss $K$-theory. A main
ingredient is the model structure on the category of \textsc{abel}ian group
objects in $\mathcal{R} (X)$ established by \textsc{Dwyer\/} and \textsc{Kan}
\cite{DwyerKan-SimplCoeff}, and generalised here to cover the equivariant
version $\mathcal{R} (W,G)$ as well (Theorem~\ref{thm:model-struct-ab}); this
allows us to characterise weak equivalences in the category $\mathcal{R}^\ab
(W,G)$ intrinsically instead of using the pullback procedure suggested in
\cite[\S2.3]{W-spaces}.

Another advantage of the model category setup is that we can give a conceptual
treatment of finiteness conditions (\S\ref{sec:finite-objects}).  An object
$Y$ of a model category~$\mathcal{C}$ is called {\it categorically finite\/}
if the representable functor $\mathcal{C} (Y, \,\cdot\,)$ commutes with
filtered colimits. For $\mathcal{R}(X)$ this recovers the usual notion of a
finite object as used in~\cite[\S2.1]{W-spaces}.

\medbreak
\noindent {\bf Theorem~\ref{thm:cof-finite}.} {\it Suppose $\mathcal{C}$ is a
 cofibrantly generated model category where all generating cofibrations have
 categorically finite domain and codomain. Then an object is cofibrant and
 categorically finite if and only if it is a retract of a finite $CW$ object.}
\medbreak

Finally, we prove that the new and old description of the linearisation map
are equivalent:

\medbreak
\noindent {\bf Theorem~\ref{thm:main}.} {\it Suppose $X$ is a pointed
  simplicial set with loop group~$G$. There is a diagram of categories,
  commuting up to natural isomorphism,
\begin{diagram}
  \mathcal{R}^\ab (X)_{\mathrm{cf}} & \rTo & \mathcal{R}^\ab (*,G)_{\mathrm{cf}} \\
  \uTo<{\tilde \bZ_X} && \uTo>{\tilde \bZ_*} \\
  \hbox to 0 pt {\hss \( \mathcal{R}_{\mathrm{f}} (X) = \) }
  \mathcal{R} (X)_{\mathrm{cf}} & \rTo & \mathcal{R} (*,G)_{\mathrm{cf}}
\end{diagram}
such that the horizontal arrows induce equivalences of $K$-theory spaces, and
the right vertical functor induces the usual linearisation map. The $K$-theory
space of $\mathcal{R} (X)_{\mathrm{cf}}$ is~$A(X)$. (The subscript
``$\mathrm{cf}$'' refers to the subcategories of cofibrant, categorically
finite objects.)}
\medbreak

Throughout the paper we will freely use the language of model categories
\cite{Quillen-HA} \cite{Hirschhorn-ModelCats} and of categories with
cofibrations and weak equivalences \cite{W-spaces}.

\section{Retractive Spaces, Presheaves and Linearisation}

\subsection{The equivariant simplex category and presheaves of sets}

We start by generalising the usual simplex category of a simplicial set to the
equivariant setting.

\begin{definition}
\label{def:simplexcat}
  Let~$G$ denote a simplicial monoid acting from the right on a simplicial
  set~$W$. The {\it equivariant simplex category\/} of~$W$, denoted
  $\simpg(W)$, is the following category: Objects are the simplices of~$W$,
  formally pairs $(n, w)$ with~$n$ a natural number and $w \in W_n$; a
  morphism $ (m, v) \rTo (n, w) $ is a pair~$ (\alpha, g) $ where~$\alpha
  \colon [m] \rTo\relax [n] $ is a morphism in~$\Delta$ and~$g$ is an element
  of~$ G_m $ such that $ v = W (\alpha)(w) \cdot g $. Composition of morphisms
  is defined as $(\beta, h) \circ (\alpha, g) := (\beta \circ \alpha, G
  (\alpha)(h) \cdot g)$.
\end{definition}

It is clear from the definition that the category~$ \mathrm{\rm Simp}_G (W) $ depends
covariantly on~$W$ and~$G$; a semi-equivariant map $ (W, G) \rTo (V, H) $ induces a functor
$ \mathrm {Simp}_G (W) \rTo \mathrm {Simp}_H (V) $.

\medbreak

There is a more conceptual description of the simplex category.  Let $H$ be a
(discrete) monoid acting from the right on a set~$S$. The {\it transport
  category\/} of~$S$ is the category $\mathrm{Tr}_H (S)$ with objects the
elements of~$S$, and a morphism $ s \rTo t $ is an element $ h \in H $ such
that $ s = t \cdot h $.  Composition of morphisms is given by multiplication
in~$H$.

The transport category is functorial in~$S$ and~$H$. In particular, a
$G$-equivariant simplicial set~$W$ gives rise to a simplicial category
$$\mathrm{Tr} \colon \Delta^\op \rTo \mathrm{Cat},\ [n] \mapsto
\mathrm{Tr}_{G_n} (W_n)\ .$$

If $\mathcal{C}$ is a category, we associate to a functor $F \colon
\mathcal{C}^\op \rTo \mathrm{Cat}$ into the category of (small) categories a
new category, the \textsc{Grothendieck\/} construction $\mathrm{Gr}(F)$
of~$F$. Objects are the pairs $(c,x)$ with $c$ an object of~$\mathcal{C}$ and
$x$ an object of~$F(c)$. A morphism $(c,x) \rTo (d,y)$ consists of a morphism
$f \colon c \rTo d$ in~$\mathcal{C}$ and a morphism $\alpha \colon x \rTo
F(f)(y)$ in~$F(c)$ with composition given by the formula $(f, \alpha) \circ
(g, \beta) := (f \circ g, F(g)(\alpha) \circ \beta)$. The category
$\mathrm{Gr}(F)$ is called the lax colimit $F \int \mathcal{C}^\op$ in
\cite[\S3]{Th:FirstQuad}.

\begin{lemma}
  \label{lem:simpgw-is-GrOr}
  The category $\simpg(W)$ is the \textsc{Grothendieck\/} construction of the
  functor $\mathrm{Tr} \colon \Delta^\op \rTo \mathrm{Cat}, [n] \mapsto
  \mathrm{Tr}_{G_n} (W_n)$.  \qed
\end{lemma}

\begin{proposition}
  \label{FunIsRofX}
  The category ${\cal R} (W, G)$ is the category of set-valued pre\-sheaves on
  $\simpg(W)$, \ie, $\mathcal{R}(W,G)$ is equivalent to $\fun (\simpg
  (W)^\op,\, \mathrm {\rm Set}_* )$. The functor corresponding to $(Y,r,s) \in
  \mathcal{R}(W,G)$ sends $(n,w) \in \simpg(W)$ to the fibre $r^{-1} (w)$ over
  $w \in W_n$. In particular, the category~$ {\cal R} (W, G) $ is complete and
  cocomplete; limits and colimits are computed fibrewise.
\end{proposition}

\begin{proof}
  An object~$ Y = (Y, r, s) \in {\cal R} (W, G) $ determines a functor
  $$ \simpg (W)^\op \rTo \mathrm{Set}_* $$
  by sending $ (n, w) $ to the set $ r^{-1} (w) $ with basepoint given by
  $s (w)$.
  
  Conversely, starting from a functor~$F$, we define the associated simplicial
  set~$Z$, in degree~$n$, as the disjoint union of sets $ \amalg_{w \in W_n} F
  (n, w) $. From the definition of~$ \simpg (W) $ one can check that~$Z$ is
  actually a $G$-equivariant simplicial set. The $G$-action is specified by
  the morphisms of type $(\mathrm{id}, g)$; explicitly, an element $y \in
  F(n,w)$ of the $w$-summand of the disjoint union is mapped to the element
  $F(1,g) (y)$ of the $wg$-summand. The simplicial structure is similarly
  encoded by morphisms of the type $(\alpha, 1)$.  The constant functor
  sending everything into the one-point-set has~$W$ as associated simplicial
  set. Thus the canonical maps to and from the one-point-set induce
  equivariant maps to and from~$W$, thereby making~$Z$ into an object
  of~${\cal R} (W, G) $.
  
  A calculation shows that these assignments are actually inverse to each
  other. We omit the details.
\end{proof}

\medbreak

  In a similar way one can show that the category of $G$-equi\-variant
  simplicial sets over~$W$ (with no specified section) is equivalent to the
  category $ \fun (\simpg (W)^\op,\, \mathrm{Set}) $.

\begin{definition}
  \label{def:ncell}
  An {\it equivariant $n$-simplex\/} or an {\it $n$-cell} in ${\cal R} (W, G)$
  is an object of $\mathcal{R}(W,G)$ isomorphic to
  $$ \Delta[n,w] = \chi_w + \mathrm{id }\colon (\Delta^n \times G) \amalg W \rTo W $$
  where $\chi_w \colon \Delta^n \times G \rTo W$ is obtained from the
  characteristic map of the simplex \hbox{$w \in W_n$} by forcing equivariance.  
  Structural section is the inclusion of~$W$ as the
  second summand. The boundary is the restriction of~$\Delta[n,w]$ to the subspace
  $(\partial \Delta^n \times G) \amalg W$, denoted~$\partial
  \Delta[n,w]$. Similarly, we can define the horns $\Lambda_i [n,w]$.
\end{definition}

\begin{proposition}
\label{SimplexRep}
The equivalence of Proposition~\ref{FunIsRofX} identifies the $n$-cell
$\Delta[n,w] \in {\cal R} (W, G)$ with the representable functor $\simpg(W)\ 
\bigl(\,\cdot\,,(n,w)\bigr)_+ $.  \qed
\end{proposition}

\medbreak

There is another relationship between an object $ Y \in {\cal R} (W, G) $ and
its associated functor $ F \colon \simpg (W)^\op \rTo \mathrm{Set}_* $ which
we describe next. Let~$H$ be a discrete group acting on a set~$S$.  We define
the {\it homotopy orbit space\/} of~$S$, denoted~$ S_{hH} $, as the nerve of
the transport category $\mathrm{\rm Tr}_H S$ of~$S$.  This definition is
functorial in~$S$ and~$H$.  If~$G$ is a simplicial group acting on a
simplicial set~$Y$, we define the {\it homotopy orbit space of~$Y$\/},
denoted~$ Y_{hG} $, as the diagonal of the bisimplicial set $ [n] \mapsto
(Y_n)_{hG_n} $. If the action of $G$ on $Y$ is free, there is a weak homotopy
equivalence $Y/G \simeq Y_{hG}$.

\begin{theorem}
  \label{ThmHocolim}
  Let~$F$ denote the functor corresponding to~$ Y \in {\cal R} (W, G) $,
  considered as a functor into the category of sets (not pointed sets). Then
  $\mathrm {hocolim}\, F \simeq Y_{hG} $.
\end{theorem}

\begin{proof}
  First note that a set $S$ determines a discrete category; objects are the
  elements of the set, the only morphisms are identities. Then $NS = S$ as
  simplicial sets (here and in what follows, $N$ denotes the nerve of a
  category).  Moreover, we can consider any set-valued functor as a functor
  into the category of (small) categories and thus speak of its
  \textsc{Grothendieck\/} construction.

  By \cite{Th:FirstQuad} (Theorem~3.19 and preceding remarks) there are weak
  equivalences
  $$ \mathrm{hocolim} (F) = \mathrm{hocolim} NF \simeq N \mathrm{Gr} (F) \ .$$
  Direct calculation shows $\mathrm{Gr}(F) \cong \simpg(Y)$. By
  Lemma~\ref{lem:simpgw-is-GrOr} the latter category
  is the \textsc{Grothendieck\/} construction of the functor
  $$ \mathrm{\rm Tr} \colon \Delta^\op \rTo Cat,\quad
  [n] \mapsto \mathrm{\rm Tr}_{G_n} Y_n $$
  which is the simplicial category of transport categories determined by~$Y$.
  Thus, using \cite[Theorem~3.19]{Th:FirstQuad} again,
  $$N \mathrm{Gr}(F) \cong N \mathrm{Gr} (\mathrm{Tr}) \simeq
  \mathrm{hocolim}_{\Delta^\op} (N \mathrm{Tr}) \ .$$
  But the functor $N\mathrm{Tr} \colon \Delta^\op \rTo \mathrm{sSet}$
  determines a bisimplicial set, and its homotopy colimit is weakly equivalent
  to its diagonal \cite[XII.4.3]{BK:Monster} which, by definition, is the
  homotopy orbit space $Y_{hG}$.
\end{proof}

\subsection{Presheaves of \textsc{abel}ian groups and linearisation}
\label{sec:presh-abel-groups}

The category $\mathcal{R}(W,G)$ is a category with products, given by pull
back over~$W$ using structural retractions, and $(W,\mathrm{id}, \mathrm{id})$
is a terminal object.  Hence it makes sense to speak of the category ${\cal
  R}^\ab (W, G)$ of \textsc{abel}ian group objects in~${\cal R} (W, G)$. For
an object $(Y,r,s) \in {\cal R}^\ab (W, G)$ the fibres of~$r$ are
\textsc{abel}ian groups, and ${\cal R}^\ab (W, G)$ is a category of presheaves
on~$W$:

\begin{proposition}
  \label{FunIsAb}
  The category $ {\cal R}^\ab (W, G) $ is equivalent to the functor
  category $ \fun(\simpg (W)^\op,\ \mathrm {\rm Ab} )$.
\end{proposition}

\begin{proof}
  Using Proposition~\ref{FunIsRofX} we have a chain of equivalences of
  categories $\mathcal{R}^\ab (W,G) \cong \fun (\simpg (W)^\op,\ \mathrm
  {Set}_* )^\ab \cong \fun (\simpg (W)^\op,\ \mathrm {Set}_*^\ab) \cong
  \fun(\simpg (W)^\op,\ \mathrm {Ab}) $.
\end{proof}

\begin{definition}
  The {\it linearisation functor\/} $\tilde\bZ_W \colon \mathcal{R}(W,G) \rTo
  \mathcal{R}^\ab (W,G)$ is the left adjoint of the forgetful functor
  $\mathcal{R}^\ab(W,G) \rTo \mathcal{R}(W,G)$.
\end{definition}

\noindent For $W = *$ this recovers \textsc{Waldhausen}'s definition of the
linearisation map \cite[p.~398]{W-spaces}. In general, the functor
$\tilde\bZ_W$ is given by ``applying $\tilde \bZ[\,\cdot\,]$ to the fibres
of~$r$''. Here $\tilde \bZ[\,\cdot\,]$ denotes the reduced free
\textsc{abel}ian group functor given by $\tilde \bZ[S] := \bZ[S]/\bZ[*]$ for a
pointed set~$S$. In the presheaf language, linearisation is particularly easy
to describe: The functor $\tilde\bZ_W$ sends $F \colon \simpg (W)^\op \rTo
\mathrm{Set}_*$ to the functor
$$\tilde\bZ_W(F) := \tilde \bZ[\,\cdot\,] \circ F \colon 
\simpg (W)^\op \rTo \mathrm{Ab},\ (n,w) \mapsto \tilde\bZ [F(n,w)] \ . $$

\medbreak

\begin{definition}
\label{def:cell-lin}
An {\it $n$-simplex\/} or an {\it $n$-cell\/} in the category $\mathcal{R}^\ab
(W,G)$ is an object isomorphic to $\Delta^\ab [n,w] := \tilde\bZ_W (\Delta [n,w])$
where $\Delta [n,w]$  is an $n$-simplex in $\mathcal{R}(W,G)$. The
boundary of $\Delta^\ab [n,w]$ is $\partial\Delta^\ab [n,w] := \tilde\bZ_W
(\partial\Delta [n,w])$. Similarly, we can define the horns $\Lambda_i^\ab
[n,w] := \tilde\bZ_W(\Lambda_i [n,w])$.
\end{definition}

\begin{proposition}
  \label{SimplexRepAb}
  The equivalence of Proposition~\ref{FunIsRofX} identifies the $n$-cell
  $\Delta^\ab[n,w]$ with the the representable functor $\tilde\bZ[\simpg(W)\ 
  \bigl(\,\cdot\,,(n,w)\bigr)_+] $.  \qed
\end{proposition}

\subsection{Functors}
\label{sec:assembly}

From now on, we suppose that $X$~is a connected pointed simplicial set, $G$ is
the loop group \cite{Kan-LoopGroup} of~$X$, and $\xi \colon P \rTo X$ is a
universal $G$-bundle, \ie, $P$~is a weakly contractible free $G$-equivariant
simplicial set, the map $\xi$ is constant on $G$-orbits, and $\xi$ induces an
isomorphism $X \cong P/G\,$. An explicit construction is given in
\cite[Lemma~9.3]{Kan-LoopGroup}.

\medbreak

We proceed to define the following diagram of categories:
\begin{diagram}[l>=3em,eqno=(\dagger)]
{\cal R}^\ab (X) & \rTo^{\Xi^*} & {\cal R}^\ab (P, G) &
        \rTo^C & {\cal R}^\ab (*, G)
\hbox to 0 pt {\(\ \cong \mathcal{M}(\bZ[G])\)\hss} \\
\uTo<{\tilde\bZ_X} && \uTo<{\tilde\bZ_P} && \uTo>{\tilde\bZ_*} \\
{\cal R} (X) & \rTo^{\overline{\Xi}^*} & {\cal R}(P, G) & \rTo^{\overline{C}} & {\cal R}
(*, G)
\hbox to 0 pt {\(\ \cong G\hbox{-}\mathrm{sSet}_*\)\hss} \\
\end{diagram}
This diagram commutes up to natural isomorphism. Note that $\mathcal{R}(*,G)$
is the category of pointed $G$-equivariant simplicial sets, and that ${\cal
  R}^\ab (*, G)$ is equivalent to the category $\mathcal{M}(\bZ[G])$ of
simplicial right $\bZ[G]$-modules.---The vertical arrows denote the
linearisation functors as defined in \S\ref{sec:presh-abel-groups}.

\medbreak

The bundle projection $\xi \colon P \rTo X$ determines a functor
$$ \Xi \colon \simpg (P)^\op \rTo \simp(X)^\op \ .$$
Since $G$ acts freely on~$P$, this functor is actually an equivalence of categories.
By pre-composition $\Xi$ determines an equivalence of categories
$$ \overline{\Xi}^* \colon \fun (\simp(X)^\op,\ \mathrm{Set}_*)
   \rTo \mathrm{Fun}\ (\simpg (P)^\op,\ \mathrm{Set}_*) \ . $$
In the language of retractive spaces, $\overline{\Xi}^*$  is given by
$$ {\cal R} (X) \rTo\relax {\cal R} (P, G), \quad Y \mapsto Y \times_X P \ .$$
The same description defines an equivalence of ${\cal
  R}^\ab (X) $ and ${\cal R}^\ab (P, G)$, denoted~$\Xi^*$.
By construction there is a natural isomorphism $\Xi^* \circ \tilde\bZ_X \cong
  \tilde\bZ_P \circ \Xi^*$.

\medbreak

The map $P \rTo *$ induces a functor $\kappa \colon \simpg (P)^\op \rTo \simpg
(*)^\op$. Both $C$ and $\overline{C}$ are defined, using the presheaf
language, by left \textsc{Kan} extension along~$\kappa$.  On the level of
retractive spaces, $\overline{C} \colon \mathcal{R}(P,G) \rTo\relax
\mathcal{R}(*,G)$ is given by collapsing $P$, sending $Y$ to $ Y/P$. For
\textsc{abel}ian group objects the collapsing functor~$C$ involves summation
over the fibres of the structural retraction. Explicitly, $C(Y)_n =
\bigoplus_{p \in P_n} r^{-1} (p)$ for an object $(Y,r,s) \in \mathcal{R}^\ab
(P,G)$. (An explicit description for $\overline{C}$ is given by the same
expression if ``$\bigoplus$'' is interpreted as a one-point union of pointed
sets.)

\section{Model Structures and Algebraic $K$-Theory}

\subsection{Categorically finite objects and finite $CW$ objects}
\label{sec:finite-objects}

An object $Y$ of a cocomplete category $\mathcal{C}$ is called {\it
  categorically finite\/} if the representable functor
$\mathcal{C}(Y,\,\cdot\,)$ commutes with filtered colimits. We want to compare
this categorical finiteness notion with the more geometric notion of a finite
$CW$ object which is defined whenever we have a notion of ``cells'' in the
category~$\mathcal{C}$.

\begin{theorem}
  \label{thm:cof-finite}
  Let $\mathcal{C}$ be a model category. Suppose there is a set $I$ of
  generating cofibrations which have categorically finite domains.  Let $Y \in
  \mathcal{C}$ be an object. If~$\,Y$ is categorically finite and cofibrant,
  then $Y$ is a retract of an object $Z$ which is a finite $CW$ object in the
  following sense: There is a finite filtration
  $$* = Z_0 \rTo Z_1 \rTo \ldots \rTo Z_n = Z$$
  where each map $Z_k \rTo Z_{k+1}$ is a pushout of a single map
  in~$I$, and ``$*$'' denotes the initial object of~$\mathcal{C}$.
  
  Conversely, if $Y$ is a retract of a finite $CW$ object~$Z$, and if the maps
  in~$I$ have categorically finite domain and codomain, then $Y$ is
  categorically finite and cofibrant.
\end{theorem}

\begin{proof}
  If $Y$ is a retract of a finite $CW$ object then $Y$ is certainly cofibrant.
  If all the maps in~$I$ have categorically finite domain and codomain, then
  since $Z$ is a finite colimit of categorically finite objects $Z$ is
  categorically finite, hence so is its retract~$Y$.
  
  To prove the first statement of the theorem, assume that $Y$ is
  categorically finite and cofibrant. Since $Y$ is cofibrant we know by
  \textsc{Quillen}'s small object argument \cite[\S{}II.3]{Quillen-HA},
  \cite[\S10.5.16]{Hirschhorn-ModelCats} that
  $Y$ is a retract of an object~$C$ which comes equipped with a filtration
  $$* = C_0 \rTo C_1 \rTo C_2 \rTo \ldots$$
  with $\mathrm{colim}\, C_i = C$
  such that each map $C_k \rTo C_{k+1}$ is a pushout of a (possibly infinite)
  coproduct of maps in~$I$. In more details, the
  small object argument yields a factorisation of the map $* \rTo Y$ as a
  cofibration $* \rTo C$ and an acyclic fibration $C
  \rTo Y$. Since $Y$ is cofibrant, we can find the dotted lift in the solid
  square diagram below exhibiting $Y$ as a retract of~$C$.
  \begin{diagram}
    * & \rTo & C \\ \dTo & \ruDotsto>i & \dTo \\ Y & \rTo_= & Y
  \end{diagram}
  
  By construction, the following square is a pushout diagram for all $k
  \geq 0$:
  \begin{diagram}
    \coprod_{\lambda \in \Lambda_k} A^k_\lambda & \rTo[l>=4em]^{\amalg f^k_\lambda} &
    \coprod_{\lambda \in \Lambda_k} B^k_\lambda \\
    \dTo<{\sum g^k_\lambda} && \dTo \\
    C_k & \rTo & C_{k+1} \\
  \end{diagram}
  Here $\Lambda_k$ is an index set, $f^k_\lambda \colon A^k_\lambda \rTo
  B^k_\lambda$ is a map in $I$, and $g^k_\lambda \colon A^k_\lambda \rTo C_k$ is the
  attaching map. For a subset $M \subseteq \Lambda_k$ define $C_k (M)$ by a
  similar pushout with indices $\lambda$ ranging through $M$ only.
  
  Since $Y$ is categorically finite, the section $i \colon Y \rTo C$ factors
  through a finite stage $C_{k+1}$, $k \geq -1$, of the filtration of~$C$. Now
  $C_{k+1}$ is obtained from $C_k$ by attaching a coproduct of cells, indexed
  by $\Lambda_k$. Hence $C_{k+1}$ is the filtered colimit of the objects
  $C_k(M)$ where $M$ varies over the (directed) poset of finite subsets of
  $\Lambda_k$. Thus there exists a finite subset $M_k \subseteq \Lambda_k$
  and a factorisation of $i$ as $Y \rTo C_k (M_k) \rTo C$.
  
  Since all maps in~$I$ have categorically finite domain by hypothesis, the
  object $\amalg_{\lambda \in M_k} A^k_\lambda$ is categorically finite. By a
  similar argument as above, the attaching map $h \colon \amalg_{\lambda \in
    M_k} A^k_\lambda \rTo C_k$ factors through $C_{k-1} (M^\prime)$ for some
  finite subset $M^\prime \subseteq \Lambda_{k-1}$. Thus $h$ factors
  compatibly through the objects $C_{k-1} (M)$ where $M$ is a finite subset of
  $\Lambda_{k-1}$ containing $M^\prime$. Consequently, $C_{k-1} (M) (M_k)$ is
  defined for such~$M$, and $C_k (M_k)$ is the filtered colimit of the
  $C_{k-1}(M)(M_k)$.  Thus there exists a finite subset $M_{k-1} \subseteq
  \Lambda_{k-1}$ containing $M^\prime$ and a factorisation of $i$ as $Y \rTo
  C_{k-1} (M_{k-1})(M_k) \rTo C$.
  
  Iterating this argument shows that $Y$ is a retract of an object that can be
  obtained from the initial object $*$ by attaching cells indexed by the
  finite set $\coprod_{i=0}^k M_i$, hence $Y$ is a retract of a finite $CW$
  object as claimed.
\end{proof}

\subsection{Model structure and simplicial structure of $\mathcal{R}(W,G)$}
\label{sec:model-nonlin}

For the rest of the paper, we will use the notational conventions of
\S\ref{sec:assembly}. That is, $X$ is a connected pointed simplicial set with
loop group~$G$, and $P \rTo X$ is a universal $G$-bundle.

\medbreak

The category $\mathcal{R}(W,G)$ admits the structure of a simplicial model
category where a map $f \colon (Y,r,s) \rTo (Y^\prime, r^\prime, s^\prime)$ is
a weak equivalence ({\it resp.}, a fibration) if and only if the underlying
map $Y \rTo Y^\prime$ of simplicial sets is a weak equivalence ({\it resp.}, a
\textsc{Kan} fibration). This model structure is cofibrantly generated; a set
of generating cofibrations is given by the inclusion maps
$$ \partial \Delta [n,w] \rTo \Delta [n,w] \quad\hbox{for}\ n \geq 0,
\ w \in W_n \ ,$$
and a set of generating acyclic cofibrations is given by the inclusion maps
$$ \Lambda_i[n,w] \rTo \Delta [n,w] \quad\hbox{for}\ n \geq 0,
\ w \in W_n,\ 0 \leq i \leq n \ .$$
Since these maps have categorically finite domain and codomain,
Theorem~\ref{thm:cof-finite} completely characterises categorically finite
cofibrant objects in~$\mathcal{R}(W,G)$.

Given $ (Y,r,s) \in \mathcal{R} (W,G)$ and a simplicial set $K$ the object
$(Y,r,s) \otimes K$ is given by the pushout $(Y \times K) \cup_{(W \times K)}
W$ (where the map $W \times K \rTo W$ is the projection), equipped with the
obvious structure maps to and from~$W$.  In particular, if $Y=\Delta[n,w]$ is
an equivariant $n$-cell (\ref{def:ncell}), then $(Y,r,s) \otimes K$ is the
object
$$ (\Delta^n \times K \times G) \amalg W \rTo W $$
with retraction induced by projection onto $\Delta^n$ and the
characteristic map of $w \in W_n$, and section given by inclusion into the
second summand. This shows:

\begin{lemma}
  \label{lem:tensor-finite}
  If $(Y,r,s) \in \mathcal{R}(W,G)$ is a finite colimit of cells
  (representable functors), so is $(Y,r,s) \otimes \Delta^k$ for all $k \geq
  0$.  \qed
\end{lemma}

\subsection{Simplicial structure of $\mathcal{R}^\ab (W,G)$}
\label{sec:simp-linear}

For two objects $Y, Z \in \mathcal{R}^\ab (W,G)$ we define their tensor
product $Y \otimes Z$ by taking fibrewise tensor product of \textsc{abel}ian
groups. In the presheaf language this is the functor
$ Y \otimes Z \colon \simpg (W)^\op \rTo \mathrm{Ab}, \ (n,w) \mapsto
Y(n,w) \otimes Z (n,w)$.

A simplicial set $K$ defines an object $(W \times K) \amalg W$ of
$\mathcal{R} (W,G)$ with structural retraction induced by projection
onto~$W$. Given an object $Y \in \mathcal{R}^\ab (W,G)$ we define
$$Y \otimes K := Y \otimes \tilde\bZ_W ((W \times K) \amalg W) \ ,$$
this determines the simplicial structure of $\mathcal{R}^\ab (W,G)$. There is
a natural isomorphism
$$ \tilde\bZ_W \left( (Y,r,s) \otimes K \right) \cong \tilde\bZ_W (Y,r,s) \otimes K $$
for all $(Y,r,s) \in \mathcal{R}(W,G)$ and simplicial sets~$K$. In
particular, $\Delta^\ab [n,w] \otimes \Delta^k$ is a finite colimit of cells
(representable functors). More generally, we have:

\begin{lemma}
  \label{lem:tensor-finite-lin}
  If $(Y,r,s) \in \mathcal{R}^\mathrm{ab}(W,G)$ is a finite colimit of cells
  (representable functors), so is $(Y,r,s) \otimes \Delta^k$ for all $k \geq
  0$.  \qed
\end{lemma}

\subsection{Model structure of $\mathcal{R}^\ab (W,G)$}
\label{sec:model-struct-ab}

Precomposition with the inclusion of categories $\iota \colon \simp (W) \rTo
\simpg (W)$ induces functors $\overline{V} \colon \mathcal{R} (W,G) \rTo
\mathcal{R} (W)$ and $V \colon \mathcal{R}^\ab (W,G) \rTo \mathcal{R}^\ab (W)$
which forget the $G$-action on retractive spaces. The left adjoints of these
are given by left \textsc{Kan} extension along~$\iota$. Explicitly, the left
adjoint of $\overline{V}$ is described by
$$ \overline{G}_* \colon \mathcal{R} (W) \rTo \mathcal{R} (W,G), \quad
Y \mapsto (Y \times G) \cup_{(W \times G)} W $$
(pushout along the action $W \times G \rTo W$ of $G$ on~$W$), and a similar
formula describes $G_*$, the left adjoint of~$V$. By general properties of
adjoints, linearisation commutes (up to canonical isomorphism) with $G_*$ and
$\overline{G}_*$. This implies:

\begin{lemma}
  The functors $G_*$ and $\overline{G}_*$ preserve cells, boundaries of cells
  and horns. More explicitly, $\overline{G}_* (\Delta[n,w]) \cong \Delta[n,w]$
  and $G_* (\Delta^\ab[n,w]) \cong \Delta^\ab[n,w]$, and similarly for
  boundaries and horns. \qed
\end{lemma}

Given an object $Y \in \mathrm{Fun} (\simpg (W)^\op, \mathrm{Ab}) \cong
\mathcal{R}^\ab (W,G)$ we define a simplicial \textsc{abel}ian group
$\bigoplus_W Y$ which in degree~$n$ is given by $(\bigoplus_W Y)_n :=
\bigoplus Y(n,w)$, the sum ranging over all $n$-simplices of~$W$ (\ie, objects
of the form $(n,w)$ of $\simp (W)$). Given $Y$ and an object $K \in
\mathcal{R} (W)$ we define
$$ \tilde H_* (K; Y) := \pi_* \bigoplus_W \left( \tilde\bZ_W K \otimes Y \right) \ .$$
(Strictly speaking one should write $V (Y)$ instead of~$Y$ in this
formula.) The construction is functorial in $K$ and $Y$.

\begin{theorem}
\label{thm:model-struct-ab}
  The category $\mathcal{R}^\ab (W,G)$ admits a simplicial model structure
  where a map $f$ is a weak equivalence if and only if $\tilde H_* (K; f)$ is an
  isomorphism for all fibrant objects $K \in \mathcal{R}(W)$, and $f$ is a
  cofibration if and only if it has the left lifting property with respect to
  all morphisms $g$ such that the underlying map of~$g$ is an acyclic
  fibration of simplicial sets.
\end{theorem}

For $W = *$ this recovers the usual model structure for simplicial
$\bZ[G]$-modules \cite[\S{}II.6]{Quillen-HA}.---The Theorem equips all the
categories in the upper row of~$(\dagger)$ with model structures.

\begin{proof}
  If~$\,G$ is the trivial group this is the main result of
  \cite[\S5]{DwyerKan-SimplCoeff}. For the general case, we can reformulate
  the definitions of weak equivalences and cofibrations. Namely, a map $f$ is
  a weak equivalence if and only if $V(f)$ is a weak equivalence in
  $\mathcal{R}^\ab (W)$, and $f$ is a cofibration if and only if it has the
  left lifting property with respect to all morphisms $g$ such that $V(g)$ is
  an acyclic fibration in $\mathcal{R}^\ab (W)$. Using this, the proof follows
  the pattern of \cite[5.7]{DwyerKan-SimplCoeff} with virtually no
  changes. For the \textsc{Bousfield} argument in
  \cite[Proposition~5.8]{DwyerKan-SimplCoeff} the cardinal $c$ has to be at
  least as large as the cardinality of $G$ and the cardinality of the total
  space of the path fibration in Lemma~4.8 of~\cite{DwyerKan-SimplCoeff}.
\end{proof}

\begin{corollary}
  The set of inclusions
  $$\partial \Delta^\ab [n,w] \rTo \Delta^\ab [n,w]\ \in
  \mathcal{R}^\ab (W,G)$$
  is a set of generating cofibrations for the model
  structure of Theorem~\ref{thm:model-struct-ab}.
\qed
\end{corollary}

Since the maps of the Corollary have categorically finite domain and co\-domain,
Theorem~\ref{thm:cof-finite} completely characterises categorically finite
cofibrant objects in~$\mathcal{R}^\ab(W,G)$.

\medbreak

We return to the diagram~$(\dagger)$. Recall the notational conventions from
the beginning of \S\ref{sec:assembly}.

\begin{proposition}
  The functor $\Xi^* \colon \mathcal{R}^\ab (X) \rTo \mathcal{R}^\ab (P,G)$
  preserves and detects weak equivalences. In particular, it is the left
  adjoint of a \textsc{Quillen} equivalence.
\end{proposition}

\begin{proof}
  Let $f \colon Y \rTo Z$ be a morphism in $\mathcal{R}^\ab (X)$. By
  \cite[Proposition~4.8]{DwyerKan-SimplCoeff}
  the map $f$ is a weak equivalence if and only if the map of
  simplicial \textsc{abel}ian groups
  $$ \bigoplus_X \left( \tilde\bZ_X (P \amalg X) \otimes Y \right) \rTo
  \bigoplus_X \left( \tilde\bZ_X (P \amalg X) \otimes Z \right) $$
  is a weak homotopy equivalence (we consider $P \amalg X$ as an object of
  $\mathcal{R}(X)$ in the obvious way). But this map can also be described as the
  map
  $$ \bigoplus_P (P \times_X Y) \rTo \bigoplus_P (P \times_X Z) $$
  induced by $\Xi^*(f)$. The domain of this map can be rewritten as
  $$ \bigoplus_P \left( \tilde\bZ_P (P \amalg P) \otimes (P \times_X Y)
  \right) = \bigoplus_P \left( \tilde\bZ_P (P \amalg P) \otimes \Xi^*(Y)
  \right) \ ,$$
  and similarly for the codomain. By \cite[Proposition~4.8]{DwyerKan-SimplCoeff}
  again, applied to the path fibration $P \rTo^= P$, this map is a weak
  homotopy equivalence if and only if $\Xi^* (f)$ is a weak equivalence
  
  The rest of the proposition follows easily since $\Xi^*$ is an equivalence
  of categories.
\end{proof}

\medbreak

Let $T$ denote the right adjoint of~$C$. It is given by
$$ T \colon \mathcal{R}^\ab (*,G) \rTo \mathcal{R}^\ab (P,G),\quad M \mapsto M
\times P$$
with structural section given by $P = \{0\} \times P \rTo M
\times P$, and retraction given by projection onto~$P$. In the presheaf
language, $T$ is given by precomposition with $\kappa \colon \simpg(P)^\op \rTo
\simpg(*)^\op$.---The functor $T$ preserves all colimits (because $T$ itself
has a right adjoint, sending $(Y,r,s)$ to the space of all
non-equivariant sections of~$r$).

\begin{proposition}
  \label{prop:adjunction-C-T}
  The two functors
  $$C \colon \mathcal{R}^\ab (P,G) \rTo \mathcal{R}^\ab (*,G) \quad \hbox{and} \quad
  T \colon \mathcal{R}^\ab (*,G) \rTo \mathcal{R}^\ab (P,G)$$
  preserve and detect weak equivalences. The unit and counit of the adjunction
  of $C$ and $T$ are natural weak equivalences. In particular, $C$ is the left
  adjoint of a \textsc{Quillen} equivalence.
\end{proposition}

\begin{proof}
  The functor $C$ preserves and detects weak equivalences by
  \cite[Proposition~4.8]{DwyerKan-SimplCoeff}, applied to the path fibration
  $P \rTo^= P$.
  
  The counit of the adjunction of $C$ and $T$ is given by the natural map
  $$\epsilon_M \colon C \circ T (M) \cong \bZ[P] \otimes_{\bZ} M \rTo \bZ[*]
  \otimes_{\bZ} M \cong M \quad \hbox {for \(M \in \mathcal{R}^\ab (*,G)\)}$$
  induced by $P \rTo *$. Since $P$ is weakly contractible, the map $f \colon
  \bZ [P] \rTo \bZ[*]$ is a weak equivalence. But both domain and codomain are
  free simplicial \textsc{abel}ian groups, thus the map is in fact a homotopy
  equivalence in the strong sense: A choice of a basepoint in $P$ determines a
  map $g \colon \bZ[*] \rTo \bZ[P]$ with $f \circ g = \mathrm{id}_{\bZ[*]}$,
  and there exists a homotopy $H \colon \bZ [\Delta^1] \otimes_\bZ \bZ [P]
  \rTo \bZ [P]$ from $\mathrm{id}_{\bZ[P]}$ to~$g \circ f$. This implies that
  $\epsilon_M = f \otimes \mathrm{id}_M$ is a homotopy equivalence of
  simplicial \textsc{abel}ian groups with homotopy inverse $g \otimes
  \mathrm{id}_M$ and homotopy $H \otimes \mathrm{id}_M$ from
  $\mathrm{id}_{\bZ[P] \otimes M}$ to $(g \otimes \mathrm{id}_M) \circ (f
  \otimes \mathrm{id}_M)$. In particular, $\epsilon_M$ induces isomorphisms on
  homotopy groups and hence is a weak equivalence in $\mathcal{R}^\ab (*,G)$.
  Moreover, since $C$ preserves and detects weak equivalences this implies that $T$
  preserves and detects weak equivalences as well.
  
  By adjointness, the composite $C \rTo^{C(\eta)} C \circ T \circ C
  \rTo^{\epsilon T} C $ is the identity, with $\eta$ the unit of the
  adjunction. Since $\epsilon_{T(M)}$ is a weak equivalence by the above,
  so is $C(\eta_M)$, hence $\eta_M$ is a weak equivalence.
\end{proof}

\subsection{Algebraic $K$-theory}

\begin{lemma}
\label{lem:preserve-finiteness}
All the functors in the diagram $(\dagger)$ preserve cofibrant categorically finite
objects.
\end{lemma}

\begin{proof}
  This follows immediately from Theorem~\ref{thm:cof-finite} together
  with the observation that all functors in~$(\dagger)$ preserve representable
  functors and colimits.
\end{proof}

\medbreak

An object $Y$ of a model category $\mathcal{C}$ is called {\it categorically
  finite up to homotopy\/} if there is a chain of weak equivalences connecting
$Y$ and a cofibrant categorically finite object. The full subcategories of
cofibrant objects which are categorically finite and categorically
finite up to homotopy, respectively, are denoted by $\mathcal{C}_\mathrm{cf}$ and
$\mathcal{C}_\mathrm{hcf}$.

\begin{proposition}
\label{prop:k-structure}
  \begin{enumerate}
  \item The functors in the diagram $(\dagger)$ induce,
    by restriction, the following diagram of categories:
    \begin{diagram}[l>=3em,eqno=(*_{\mathrm{\rm hcf}})]
      {\cal R}^\ab (X)_{\mathrm{\rm hcf}} & \rTo^{\Xi^*} & {\cal
        R}^\ab (P, G)_{\mathrm{\rm hcf}} &
      \rTo^C & {\cal R}^\ab (*, G)_{\mathrm{\rm hcf}} \\
      \uTo<{\tilde\bZ_X} && \uTo<{\tilde\bZ_P} && \uTo>{\tilde\bZ_*} \\
      {\cal R} (X)_{\mathrm{\rm hcf}} & \rTo^{\overline{\Xi}^*} & {\cal R}(P,
      G)_{\mathrm{\rm hcf}} & \rTo^{\overline{C}} & {\cal R}
      (*, G)_{\mathrm{\rm hcf}} \\
    \end{diagram}
    There is also a similar diagram, denoted $(*_{\mathrm{cf}})$, using
    categorically finite objects throughout. Both are diagrams of categories
    with cofibrations and weak equivalences, and exact functors.
  \item The inclusion of the diagram $(*_{\mathrm{cf}})$ into
    $(*_{\mathrm{hcf}})$ induces a homotopy equivalence of $K$-theory spaces
    for each entry.
  \end{enumerate}
\end{proposition}

\begin{proof}
  (1) It has been observed before that all the functors in the diagram
  $(\dagger)$ preserve categorical finiteness of cofibrant objects
  (Lemma~\ref{lem:preserve-finiteness}), weak equivalences and cofibrations,
  hence they preserve cofibrancy and categorical homotopy finiteness. All the
  categories are categories with cofibrations and weak equivalences; the only
  non-trivial thing to verify is that for a diagram $A \lTo B \rTo^i C$ of
  objects categorically finite up to homotopy with $i$ a cofibration, the pushout $A
  \cup_B C$ is categorically finite up to homotopy. This follows from
  \cite[Proposition~3.2]{Sagave-KModel}.  The hypothesis that $\,\cdot\,
  \otimes \Delta^1$ preserves finite objects is verified in
  Lemmas~\ref{lem:tensor-finite} and~\ref{lem:tensor-finite-lin}. Note that
  the Proposition applies to $\mathcal{R}^\ab (X)$ (and hence to the
  equivalent category $\mathcal{R}^\ab (P,G)$) although we do not have a set
  of generating acyclic cofibrations. There is a fibrant replacement functor
  defined by the process of ``filling (linearised) horns''
  \cite[Proposition~4.6]{DwyerKan-SimplCoeff}; this is enough for the proof of
  \cite[Proposition~3.2]{Sagave-KModel}.

  (2) This is part of \cite[Proposition~3.2]{Sagave-KModel}.
\end{proof}

\medbreak

\begin{lemma}
\label{lem:T-preserves}
  The functor $T$ (defined in \S\ref{sec:model-struct-ab}) preserves
  cofibrations, weak equivalences between cofibrant objects, and categorically
  homotopy finite cofibrant objects.
\end{lemma}

\begin{proof}
  The functor $T$ preserves all weak equivalences by
  Proposition~\ref{prop:adjunction-C-T}.

  Taking product with~$P$ defines a functor $\overline{T} \colon
  \mathcal{R}(*,G) \rTo \mathcal{R}(P,G)$ with $\tilde\bZ_P \circ \overline{T}
  = T \circ \tilde\bZ_*$ (on the level of presheaves, both $T$ and
  $\overline{T}$ are giving by precomposition with $\kappa \colon
  \simpg(P)^\op \rTo \simpg(*)^\op$).  Since $P$ has a free $G$-action, the
  functor $\overline{T}$ maps generating cofibrations to cofibrations. Since
  the generating cofibrations in $\mathcal{R}^\ab (*,G)$ are the images of the
  generating cofibrations in $\mathcal{R}(*,G)$ under the
  functor~$\tilde\bZ_*$, and since $\tilde\bZ_P$ preserves cofibrations, this
  implies that $T$ preserves cofibrations.
  
  To prove that $T$ preserves cofibrant objects which are categorically finite
  up to homotopy, it is enough to show that $T$ maps finite $CW$ objects
  (\ref{thm:cof-finite}) to objects categorically finite up to homotopy. We
  first consider the effect of $T$ on a single cell $M:= \tilde \bZ_*[(\Delta^n
  \times G)_+] \in \mathcal{R}^\ab ( *,G)$.  Choose a simplex $p \in P_n$.
  Since $C$ preserves cells (representable functors) we have $C
  (\Delta^\ab[n,p]) \cong M$, and the counit defines a weak equivalence
  (Proposition~\ref{prop:adjunction-C-T}) 
  $$ T(M) \cong T \circ C (\Delta^\ab[n,p])
  \rTo^{\epsilon} \Delta^\ab[n,p] $$
  which proves that $T(M)$ is categorically finite up to homotopy.
  
  Given a pushout diagram $A \lTo B \rTo^i C$ of cofibrant objects in
  $\mathcal{R}^\ab (*,G)$ with $i$ a cofibration, then if $T(A)$, $T(B)$
  and~$T(C)$ are categorically finite up to homotopy, so is $T (A \cup_B C) \cong
  T(A) \cup_{T(B)} T(C)$ since $\mathcal{R}^\ab (P,G)_{\mathrm{hcf}}$ is
  closed under pushouts by Proposition~\ref{prop:k-structure}~(1). Since the
  boundary of a cell is a finite colimit of cells of lower dimensions, the
  result now follows by a double induction on $n$ and the number of cells.
\end{proof}

\begin{theorem}
  \label{thm:main}
  All horizontal functors in the diagrams $(*_{\mathrm{cf}})$ and
  $(*_{\mathrm{hcf}})$ induce equivalences of $K$-theory spaces. In
  particular, the composite $C \circ \Xi^*$ (\S\ref{sec:assembly}) induces a
  homotopy equivalence $K (\mathcal{M}(\bZ[G])_{\mathrm{cf}}) \simeq \Omega |
  \mathcal{S}_\bullet h\mathcal{R}^\ab (X)_{\mathrm{cf}}|$, and the
  linearisation functor $\tilde\bZ_X$ induces a map
  $$A(X) \rTo \Omega |\mathcal{S}_\bullet h\mathcal{R}^\ab (X)_{\mathrm{cf}}|$$
  which is identified by $C \circ \Xi^*$ and $\overline{C} \circ
  \overline{\Xi}^*$ with the usual linearisation map in algebraic $K$-theory
  \cite[p.~398]{W-spaces} defined using~$\tilde\bZ_*$.
\end{theorem}

\begin{proof}
  First observe that $\mathcal{R}(X)_{\mathrm{cf}} =
  \mathcal{R}_{\mathrm{f}}(X)$, the latter being the full subcategory of
  $\mathcal{R}(X)$ whose objects are obtained from~$X$ by attaching finitely
  many cells \cite[\S2.1]{W-spaces}. In fact, for $Y \in
  \mathcal{R}(X)_{\mathrm{cf}}$ we know by Theorem~\ref{thm:cof-finite} that
  $Y/X$ contains finitely many non-degenerate simplices, and we can
  reconstruct $Y$ by attaching these simplices to~$X$ (no retracts necessary).
  Thus $\Omega |h\mathcal{S}_\bullet \mathcal{R}(X)_{\mathrm{cf}}| =
  A(X)$.
  
  For the remaining statements, it suffices by
  Proposition~\ref{prop:k-structure}~(2) to consider the
  diagram~$(*_{\mathrm{hcf}})$. The functors $\overline{\Xi}^*$ and~$\Xi^*$
  are exact equivalences of categories, and~$\overline{C}$ induces a homotopy
  equivalence on $K$-theory spaces by the argument of
  \cite[Proposition~2.1.4]{W-spaces}: Since $P \simeq *$, the functor
  $\overline{T}$ (cf.~Proof of Lemma~\ref{lem:T-preserves}) provides a
  homotopy inverse.
 
  We are left to consider the functor~$C$.  By the previous lemma the right
  adjoint $T$ of $C$ induces a map on $K$-theory spaces. Since unit and counit
  of the adjunction are both weak equivalences by
  Proposition~\ref{prop:adjunction-C-T}, $T$ and $C$ induce mutually inverse
  homotopy equivalences on $K$-theory spaces.
\end{proof}

\subsubsection*{Acknowledgements}

The author has to thank \textsc{S.~Sagave} and \textsc{O.~Renaudin} for
helpful discussions. The referee's comments helped significantly to improve
the exposition of the paper.

\bibliographystyle{plain}

\begin{thebibliography}{1}

\bibitem{BK:Monster}
A.~K. Bousfield and D.~M. Kan.
\newblock {\em Homotopy limits, completions and localizations}.
\newblock Springer-Verlag, Berlin, 1972.
\newblock Lecture Notes in Mathematics, Vol. 304.

\bibitem{DwyerKan-SimplCoeff}
W.~G. Dwyer and D.~M. Kan.
\newblock Homology with simplicial coefficients.
\newblock In {\em Algebraic topology (Arcata, CA, 1986)}, volume 1370 of {\em
  Lecture Notes in Math.}, pages 143--149. Springer, Berlin, 1989.

\bibitem{Hirschhorn-ModelCats}
Philip~S. Hirschhorn.
\newblock {\em Model categories and their localizations}, volume~99 of {\em
  Mathematical Surveys and Monographs}.
\newblock American Mathematical Society, Providence, RI, 2003.

\bibitem{Kan-LoopGroup}
Daniel~M. Kan.
\newblock A combinatorial definition of homotopy groups.
\newblock {\em Ann. of Math. (2)}, 67:282--312, 1958.

\bibitem{Quillen-HA}
Daniel~G. Quillen.
\newblock {\em Homotopical algebra}.
\newblock Springer-Verlag, Berlin, 1967.

\bibitem{Sagave-KModel}
Steffen Sagave.
\newblock On the algebraic {$K$}-theory of model categories.
\newblock {\em J. Pure Appl. Algebra}, 190(1-3):329--340, 2004.

\bibitem{Th:FirstQuad}
Robert~W. Thomason.
\newblock First quadrant spectral sequences in algebraic {$K$}-theory via
  homotopy colimits.
\newblock {\em Comm. Algebra}, 10(15):1589--1668, 1982.

\bibitem{W-spaces}
Friedhelm Waldhausen.
\newblock Algebraic {$K$}-theory of spaces.
\newblock In Andrew Ranicki, N.~Levitt, and F.~Quinn, editors, {\em Algebraic
  and Geometric Topology}, number 1126 in Lecture Notes in Mathematics, pages
  318--419, Rutgers, 1983. Springer-Verlag.

\end{thebibliography}

\end{document}